\documentstyle[12pt]{article}

\newtheorem{theo}{Theorem}
\newtheorem{prop}[theo]{Proposition}
\newtheorem{coro}[theo]{Corollary}
\newtheorem{lemme}[theo]{Lemma}
\newtheorem{defi}[theo]{Definition}

\def\NN{\hbox{\sf I\hskip -1pt N}}
\def\pNN{\hbox{{\scriptsize \sf I}\hskip -1pt {\scriptsize \sf N}}}
\def\ZZ{\hbox{\sf Z\hskip -4pt Z}}
\def\pZZ{\hbox{{\scriptsize \sf Z}\hskip -4pt {\scriptsize \sf Z}}}
\def\QQ{\hbox{\sf I\hskip -6pt Q}}
\def\RR{\hbox{\sf I\hskip -2pt R}}

\def\card{\hbox{\rm Card}}
\def\min{\hbox{\rm min}}
\def\max{\hbox{\rm max}}
\def\inf{\hbox{\rm inf}}

\def\mvide{{\epsilon }}
\def\gcd{\hbox{\rm gcd}}
\def\scm{\hbox{\rm scm}}
\def\id{\hbox{\rm Id}}

\def\B{{\cal B}}

\def\F{{\cal F}}

\def\M{{\cal M}}
\def\O{{\cal O}}
\def\P{{\cal P}}
\def\R{{\cal R}}

\def\xb{\overline{x}}
\def\yb{\overline{y}}

\def\fleche #1{\buildrel{#1}\over\longrightarrow}

\begin{document}

\title
{
\bf Linearly recurrent subshifts have a finite number of non-periodic subshift factors
}

\author{
Fabien Durand
}


\maketitle

\vskip 1.5cm

{\bf Abstract.} A minimal subshift $(X,T)$ is linearly recurrent if there exists a constant $K$ so that for each clopen set $U$ generated by a finite word $u$ the return time to $U$, with respect to $T$, is bounded by $K|u|$. We prove that given a linearly recurrent subshift $(X,T)$ the set of its non-periodic subshift factors is finite up to isomorphism. We also give a constructive characterization of these subshifts.

\section{Introduction.}

In this paper we continue the study of linearly recurrent (LR) subshifts initiated in \cite{DHS}. A minimal subshift $(X,T)$ is linearly recurrent if there exists a constant $K$ so that for each clopen set $U$ generated by a finite word $u$ the return time to $U$, with respect to $T$, is bounded by $K|u|$. We focus our attention on their topological Cantor factors. The present work is motivated by the fact, proved in \cite{DHS}, that given a LR subshift $(X,T)$ there exists a constant $D$ such that for all chains of non-periodic subshift factors 
$$ 
(X,T)\; \fleche{ \gamma_{n-1}} \;
(Y_{n-1},T)\;  \fleche{ \gamma_{n-2}}\; 
\cdots\; \fleche{\gamma_1}(Y_1,T)\;
\fleche{\gamma_0} \; (Y_0,T)\  ,
$$
if $n\geq D$ then there exists an integer $i$ for which $\gamma_i$ is an isomorphism. A natural question arise: Do LR subshifts have a finite number of subshift factors (up to isomorphism) ? The main result of this paper answers to this question:
\begin{theo}
\label{maintheo}
Let $(X,T)$ be a linearly recurrent subshift. The set of its non-periodic subshift factors is finite up to isomorphism.
\end{theo}
The proof of this result intensively uses the notion of return words and their properties established in \cite{Du1,DHS}. In Section \ref{sectionun} we recall these properties. In Section \ref{sectiondeux} we characterize LR subshifts by means of $S$-adic subshifts, this notion was introduced in \cite{Fe}. We also give a necessary and sufficient condition for a sturmian subshift to be LR, and we prove that sturmian subshifts $(X,T)$ are Cantor prime, i.e. they do not have proper Cantor factor, unless it is well-known they are not prime. This last result was certainly known but the author did not find reference. In Section \ref{sectiontrois} we established that LR subshifts are uniquely ergodic. Section \ref{sectionquatre} is devoted to the proof of Theorem \ref{maintheo}. We take the opportunity to make some comments on the result we recalled before concerning chains of factors, we show that LR subshifts are coalescent and that if two LR subshifts are weakly isomorphic then they are isomorphic. In the last section we consider the case of substitution subshifts and we give a necessary and sufficient condition for such subshifts to have a finite number (up to isomorphism) of Cantor factors and to have a finite number (up to isomorphism) of non-periodic Cantor factors.

\section{Definitions and background.}
\label{sectionun}

\subsection{Words, sequences and morphisms.}

We call {\it alphabet} a finite set of elements called {\it letters}. Let $A$ be an alphabet, a {\it word} on $A$ is an element of the free mono\"{\i}d generated by $A$, denoted by $A^*$, i.e. a finite sequence (possibly empty) of letters. Let $x = x_0x_1 \cdots x_{n-1}$ be a word, its {\it length} is $n$ and is denoted by $|x|$. The {\it empty-word} is denoted by $\mvide$, $|\mvide| = 0$. The set of non-empty words on  $A$ is denoted by $A^+$. If $J= [i,j]$ is an interval of $\NN = \{ 0,1\cdots \}$ then $x_J$ denote the word $x_i x_{i+1}\cdots x_j$ and is called a {\it factor} of $x$. We say that $x_J$ is a prefix of $x$ when $i=0$ and a suffix when $j=n-1$. If $u$ is a factor of $x$, we call {\it
occurrence} of $u$ in $x$ every integer $i$ such that $x_{[i,i + |u| - 1]}= u$. Let $u$ and $v$ be two words, we denote by $L_u(v)$ the number of occurrences of $u$ in $v$. 

The elements of $A^{\pZZ}$ and $A^{\pNN}$ are respectively called {\it sequences} and {\it one-sided sequences}. For a sequence $x = (x_n ; n\in \ZZ) = \cdots x_{-1}.x_0 x_1\cdots $ or a one sided-sequence $x = (x_n ; n\in \NN) = x_0 x_1\cdots $ we use the notation $x_J$ and the terms ``occurrence'' and ``factor'' exactly as for a word. We set $x^+ = x_0 x_1 \cdots $ and $x^- = \cdots x_{-2} x_{-1}$. The set of factors of length $n$ of $x$ is written $L_n(x)$, and the set of factors of $x$, or {\it language} of $x$, is represented by $L(x)$. 

The sequence $x$ is {\it periodic} if it
is the infinite concatenation of a word $v$. A {\it gap} of a factor $u$ of $x$ is an integer $g$ which is the difference between two successive occurrences of $u$ in $x$. We say that $x$ is {\it uniformly recurrent} if each factor has bounded gaps.

Let $A$, $B$ and $C$ be three alphabets. A {\it morphism} $\tau$ is a map from $A$ to $B^*$. Such a map induces by concatenation a map from $A^*$ to $B^*$. If $\tau (A)$ is included in $B^+$, it induces a map from $A^{\pNN}$ to $B^{\pNN}$ and a map from $A^{\pZZ}$ to $B^{\pZZ}$ defined by $\tau (x^-.x^+) = \tau (x^-).\tau (x^+)$. All these maps are written $\tau$ also.

To a morphism $\tau$, from $A$ to $B^*$, is na\-tu\-rally associa\-ted the ma\-trix $M_{\tau} = (m_{i,j})_{i\in  B  , j \in  A  }$ where $m_{i,j}$ is the number of occurrences of $i$ in the word $\tau(j)$. To the composition of morphisms corresponds the multiplication of matrices. For example, let $\tau_1 : B \rightarrow C^*$, $\tau_2 : A \rightarrow B^*$ and $\tau_3 : A \rightarrow C^*$ be three morphisms such that $\tau_1 \tau_2 = \tau_3$, then we have the following equality: $M_{\tau_1} M_{\tau_2} = M_{\tau_3}$. In particular if $\tau $ is a morphism from $A$ to $A^{*}$ we have $M_{\tau^n} = M_{\tau}^n$ for all non-negative integers $n$.

\subsection{Dynamical systems and subshifts.}

By a {\it dynamical system} we mean a pair $(X,S)$ where $X$ is a compact metric space and $S$ a homeomorphism from $X$ onto itself. We say that it is a {\it Cantor system} if $ X $ is a Cantor space. That is, $ X $ has a countable basis of its topology which consists of closed and open sets and does not have isolated points. The system $(X,S)$ is {\it minimal} whenever $X$ and the empty set are the only $S$-invariant closed subsets of $X$. We say that a minimal system $(X,S)$ is {\it periodic} whenever $X$ is finite. We say it is $p$-periodic if $\card (X) = p$.

Let $(X,S)$ and $(Y,T)$ be two dynamical systems. We say that $(Y,T)$ is a {\it factor} of $(X,S)$ if there is a continuous and onto map $\phi : X \rightarrow Y$ such that $\phi S = T \phi$ ($\phi $ is called {\it factor map}). If $\phi$ is one-to-one we say that $\phi $ is an {\it isomorphism} and that $(X,S)$ and $(Y,T)$ are {\it isomorphic}. We say that $(Y,T)$ (resp. $\phi$) is a {\it proper factor} (resp. {\it proper factor map}) if $(Y,T)$ is not isomorphic to $(X,S)$ and $Y$ is not reduced to one point (resp. if $\phi$ is not an isomorphism and does not identify all points of $X$). 

Let $(X,S)$ be a minimal Cantor system and $U\subset X$ be a clopen set. Let $S_U : U\rightarrow U$ be the {\it induced transformation}, i.e. if $x\in U$ then
\begin{center}
$S_U (x) = S^{r_U(x)} (x) $, where $r_U (x) = \inf \{ n > 0 \ ; \  S^n (x) \in U \}$.
\end{center}
The pair $(U , S_U )$ is a minimal Cantor system, we say that $(U , S_U)$ is the {\it induced system of $(X,S)$ with respect to $U$}.

In this paper we deal with Cantor systems called {\it subshifts}. Let $A$ be an alphabet. We endow $A^{\pZZ}$ with the topology defined by the metric
$$
d(x,y) = \frac{1}{2^n} \ \ {\rm with} \ \ n = \inf \{ |k| ; x_k \not = y_k \},
$$
where $x = (x_n ; n\in \ZZ )$ and $y = (y_n ; n\in \ZZ )$ are two elements of $A^{\pZZ}$. By a {\it subshift} on $A$ we shall mean a couple $(X,T_{/X})$ where $X$ is a closed $T$-invariant ($T(X) = X$) subset of $A^{\pZZ}$ and $T$ is the {\it shift transformation} 
\begin{center}
\begin{tabular}{lllll}
$T$ & : & $A^{\pZZ}$            & $\rightarrow $ & $A^{\pZZ}$ \\
    &   & $(x_n ; n \in \ZZ )$ & $\mapsto$      & $(x_{n+1} ; n \in \ZZ )$.
\end{tabular}
\end{center}
We call language of $X$ the set $L (X) = \{ x_{[i,j]} ; x\in X , i\leq j\}$.
Let $u$ and $v$ be two words of $A^{*}$. The set 
$$
[u.v]_X = \{ x\in X ; x_{[-|u|,|v|-1]} = uv \}
$$
is called {\it cylinder}. The family of these sets is a base of the induced topology on $X$. When it will not create confusion we will write $[u.v]$ and $T$ instead of $[u.v]_{X}$ and $T_{/X}$. We set $[v] = [\mvide . v]$.

Let $x$ be a sequence (or a one-sided sequence) on $A$ and $\Omega (x)$ be the set $\{ y \in A^{\pZZ} ; y_{[i,j]} \in L(x), \forall \ [i,j]\subset \ZZ \}$. It is clear that $(\Omega (x), T)$ is a subshift. We say that $(\Omega (x), T)$ is the subshift  generated by $x$. When $x$ is a sequence we have $\Omega (x) = \overline{\{ T^n x ; n\in \ZZ \}}$. Let $(X,T)$ be a subshift on $A$, the following are equivalent :

$\imath )$ $(X,T)$ is minimal.

$\imath\imath )$ For all $x\in X$ we have $X=\Omega (x)= \Omega (x^+)$.

$\imath\imath\imath )$ For all $x\in X$ we have $L(X)=L(x)= L(x^+)$.

We also have that $(\Omega (x) , T)$ is minimal if and only if $x$ is uniformly recurrent.

\subsection{Linearly recurrent subshifts.}

Let $A$ be an alphabet and $x$ be a sequence (or one-sided sequence) on $A$. Let $u$ be a word of $L (x)$. We call {\it return word} to $u$ of $x$ each word $w$ such that $wu$ belongs to $L (x)$, $u$ is a prefix of $wu$ and $u$ has exactly 2 occurrences in $wu$ (this notion was introduced in \cite{Du1}). We denote by $\R_{x,u}$ the set of return words to $u$ of $x$. There is a more general definition of return words in \cite{DHS} (See also \cite{Du2,HZ}). 

When $x$ is a uniformly recurrent sequence on $A$ it is easy to see that for all $u\in L (x)$ the set $\R_{x,u}$ is finite. 

\begin{defi}
We say that $x$ is {\rm linearly recurrent} (LR) (with constant $K\in\NN$) if it is uniformly recurrent and if for all $u\in L(x)$ and all $w \in \R_{x,u}$ we have $|w| \leq K|u|$.
We say that a subshift $(X,T)$ is {\rm linearly recurrent} (with constant $K$) if it is minimal and contains a LR sequence (with constant $K$). 
\end{defi}

Hence a minimal subshift is LR if and only if all its elements are linearly recurrent. These systems were introduced in \cite{DHS}. For all $x,y\in X$ we have $\R_{x,u} = \R_{y,u}$ hence we set $\R_{X,u} = \R_{x,u}$.

\subsection{Some properties of return words and LR subshifts.}

Let $A$ be an alphabet and $x$ be a uniformly recurrent sequence of $A^{\pZZ}$. Let $u$ be a prefix of $x^+$. The sequence $x$ can be written naturally as a concatenation
$$
x = \cdots m_{-2}m_{-1} . m_0m_1m_2 \cdots \:\: , m_i \in \R_{x,u}, \:\: i\in \ZZ,
$$
of return words to $u$ and this decomposition is unique. We enumerate the elements of $\R_{x,u}$ in the order of their first appearance in $x^+$. This defines a bijective map
$$
\Theta_{x,u} : R_{x,u}  \rightarrow {\R}_{x,u} \subset A^{*}
$$
where $R_{x,u} = \{ 1, \cdots , \card({\R}_{x,u})\}$. The map $\Theta_{x,u}$ defines a morphism from $R_{x,u}$ to $A^*$ and the set $\Theta_{x,u} (R_{x,u} ^*)$ consists of all concatenations of return words to $u$. When it does not create confusion we will forget the ``$x$'' in the symbols $\Theta_{x,u}$, $\R_{x,u}$ and $R_{x,u}$.

The following proposition points out some properties of return words and LR subshifts which were established in \cite{Du1} and \cite{DHS}.

\begin{prop}[\cite{Du1}]
\label{derder}
Let $A$ be an alphabet, $x$ be a uniformly recurrent sequence of $A^{\pZZ}$, $u$ be a non-empty prefix of $x^+$ and $v$ be a prefix of $u$.
\begin{enumerate}
\item 
The set $\R_{x,u}$ is a code, i.e. $\Theta_{x,u} : R_{x,u}^* \rightarrow \Theta_{x,u} (R_{x,u} ^*)$ is one to one. 
\item 
Each return word to $u$ belongs to $\Theta_{x,v} (R_{x,v} ^*)$, i.e. it is a concatenation of return words to $v$.
\item
There exists a unique map $\lambda$ from $R_{x,u}$ to $R_{x,v}^*$ such that $\Theta_{x,v} \lambda = \Theta_{x,u}$.
\end{enumerate}
\end{prop}

\begin{prop}[\cite{Du1}]
\label{lemtech}
Let $x$ be a non-periodic uniformly recurrent sequence, then
$$
\min\{ | v | ; v\in \R_{x,u} , |u| = n \} \rightarrow +\infty \;\; {\rm when} \;\; n\rightarrow +\infty.
$$
\end{prop}

\begin{prop}[\cite{DHS}]
\label{proplr}
Let $(X,T)$ be an aperiodic LR subshift with constant $K$. Then:                  
\begin{enumerate}
\item
For all $n\in \NN$ each word of length $n$ has an occurrence in each word of length $(K+1)n$.

\item
The number of distinct factors of length $n$ in 
$L(X)$ is less or equal to $Kn$.  

\item 
$X$ is $(K+1)$-power free
(i.e. $u^{K+1}\in  L(X)$ if and only if
$u=\emptyset$).  
 
\item
For all $ u\in  L(X)$ and for all $w \in\R_{u}$ we have $(1/K)|u| < |w| $.   
                
\item
For all $u\in
 L(X)$, $\card ( \R_u) \leq K(K+1)^2$.
\end{enumerate}
\end{prop}

\subsection{$S$-adic subshifts.}

Let $A$ be an alphabet, $a$ be a letter of $A$, $S$ a finite set of morphisms $\sigma$ from $A(\sigma)\subset A$ to $A^{*}$ and $(\sigma_{n} : A_{n+1} \rightarrow A_{n}^{*}; n\in \NN )$ be a sequence of $S^{\pNN}$ such that $(\sigma_{0}\sigma_{1}\cdots\sigma_{n} (aa\cdots ) ; n\in\NN)$ converges in $A^{\pNN}$ to $x$. We will say that $x$ is a {\it $S$-adic sequence on $A$}  (generated by $(\sigma_i ; i\in \NN) \in S^{\pNN}$ and $a$). The subshift generated by $x$ is called $S$-{\it adic subshift} (on $A$ generated by $x$). Of course we can always suppose that, for all $n\in \NN$, each letter of $A_n$ has an occurrence in some $\sigma_n (b)$ with $b\in A_{n+1}$. We will do so.

Let $A$ be an alphabet, $S$ be a finite set of morphisms as above and $a\in A$ such that for all $\sigma\in S$ the first letter of $\sigma (a)$ is $a$. If $(\sigma_n ; n\in \NN )$ is an element of $S^{\pNN}$ such that $\lim_{n\rightarrow +\infty} | \sigma_0 \sigma_1 \cdots \sigma_n (a) | = +\infty$, then $(\sigma_0 \sigma_1 \cdots \sigma_n (aa\cdots) ; n\in \NN )$ converges in $A^{\pNN}$. In the sequel we will always suppose that $S$-adic sequences are obtained in this way.

Let $A$ be an alphabet and $x$ be a $S$-adic sequence on $A$ generated by $(\sigma_n : A_{n+1} \rightarrow A_{n}^{*}; n\in \NN )$ and $a$. If there exists an integer $s_0$ such that for all non-negative integers $r$ and all $b\in A_{r}$ and $c\in A_{r+s_0+1}$, the letter $b$ has an occurrence in $\sigma_{r+1}\sigma_{r+2}\cdots \sigma_{r+s_0} (c)$, then we say that $x$ is a {\it primitive} $S$-adic sequence (with constant $s_0$).

\subsection{Substitution subshifts and odometers.}

A {\it substitution} on the alphabet $A$ is a morphism $\sigma : A\to A^*$ satisfying

{\parindent = 2cm 
$\imath$) There exists $a\in A$ such that $a$ is the first letter of $\sigma(a)$;

$\imath\imath$) For all $b\in A$, $\lim_{n\rightarrow +\infty} |\sigma^n (b)| = +\infty$.
\par}

It is classical that $( \sigma^n (aa\cdots) ; n\in \NN )$ converges in $A^{\pNN}$ to a sequence $x$. This sequence is a {\it fixed point} of $\sigma$, i.e. $\sigma (x) = x$.

In this paper we only consider {\it primitive} substitutions, i.e. substitutions which matrices are primitive. We recall that a matrix $M$ is primitive if and only if there exists an integer $n$ such that all the coefficients of $M^n$ are positive. In this case all the fixed points of $\sigma$ are uniformly recurrent and generate the same minimal subshift, we call it the {\it substitution subshift generated by $\sigma$}. (For more details see \cite{Qu} and \cite{DHS}.)

We remark that substitution subshifts are primitive $S$-adic subshifts. It is proved in \cite{DHS} that they are LR.

From now on we will forget to mention that the substitutions we consider are primitive.

Let $(d_n ; n\in \NN )$ be a sequence of positive integers. The inverse limit of the sequence of groups $(\ZZ / d_0d_1 \cdots d_n \ZZ \ ; \ n\in \NN )$ endowed with the addition of 1 is called {\it odometer} with base $(d_n ; n\in \NN )$. In other words an odometer with base $(d_n ; n\in \NN )$ is the system $(X,S)$ where $X = \Pi_{n=0}^{+\infty} \{ 0,1,\cdots , d_n - 1 \}$ and $S : X \to X$ is defined by
\begin{center}
\begin{tabular}{lll}
$S ( d_0 -1 , d_1 - 1 , \cdots )$ & $=$ & $(0,0,\cdots )$ and \\
$S( x_0 , x_1,\cdots )$           & $=$ & $(y_0,y_1,\cdots )$ where \\
                                  &     & $y_i=0$ if $0\leq i < i_0 = \inf \{ n\geq 0 ; x_n \not = d_n - 1 \}$, \\
                                  &     & $y_{i_0} = x_{i_0} + 1$ and 
                                          $y_i = x_i$ if $i> i_0$, elsewhere.
\end{tabular}
\end{center}

If there exists an integer $n_0$ such that for all $n\geq n_0$ we have $d_n = d_{n+1}$ then we will say that $(X,S)$ is an {\it odometer with stationary base} (it is easy to check that $(X,S)$ is isomorphic to the odometer with base $( d_0 d_1\cdots d_{n_0 - 1} , d_{n_0} , d_{n_0} , \cdots )$).

\section{A subshift is LR if and only if it is a primitive $S$-adic subshift.}
\label{sectiondeux}

The main result in this section is the following.
\begin{prop}
\label{sadiqlinrec}
Let $(X,T)$ be a subshift. The subshift $(X,T)$ is a primitive $S$-adic subshift if and only if it is LR.
\end{prop}
The proof of this proposition is inspired by the technic used in \cite{DHS} to prove that substitution subshifts are LR.

\subsection{Proof of Proposition \protect{\ref{sadiqlinrec}}.}

\begin{lemme}
\label{sadiquelb}
Let $(X,T)$ be a $S$-adic subshift on $A$ generated by $(\sigma_n : A_{n+1} \rightarrow A_{n}^* ; n\in \NN ) \in S^{\pNN}$ and $a\in A$. If for all non-negative integers $r$ there exists an integer $s\geq r$ such that for all letters $b\in A_r$ and $c\in A_{s+1}$ the letter $b$ has an occurrence in $\sigma_{r}\sigma_{r+1}\cdots \sigma_{s} (c)$, then the sequence $x = \lim_{n\rightarrow +\infty} \sigma_0\sigma_1 \cdots \sigma_n (a)$ is uniformly recurrent and $(X,T)$ is minimal.
\end{lemme}
{\bf Proof.} Let $u$ be a word of $L(X)$ and $x = \lim_{n\rightarrow +\infty} \sigma_0\sigma_1 \cdots \sigma_n (a)$. It suffices to prove that $u$ has bounded gaps in $x$. 

For each non-negative integer $r$ the word $\sigma_0\sigma_1 \cdots \sigma_r (a)$ is a prefix of $x$, hence there exists a non-negative integer $r$ such that $u$ has an occurrence in $\sigma_0\sigma_1 \cdots \sigma_r (a)$. The sequence ($\sigma_{r+1}\sigma_{r+2} \cdots \sigma_n (aa\cdots ) ; n\geq r+1)$ converges to a sequence we call $y$. 

Let $b$ be a letter having an occurrence in $y$. There exists an integer $n$ such that for each letter $c$ of $A_{n-1}$ the letter $b$ has an occurrence in $\sigma_{r+1}\sigma_{r+2} \cdots \sigma_n (c )$. Hence if $i$ and $j$ are two successive occurrences of $b$ in $y$ then
$$
|i-j|\leq 2\max \{ |\sigma_{r+1}\sigma_{r+2}\cdots \sigma_{n} (c)| ; c\in A\} = K.
$$

Now, let $i$ and $j$ be two successive occurrences of $u$ in $x$, then we have
$$
|i-j|\leq \max \{ |\sigma_{0}\sigma_1\cdots \sigma_{r} (w)| ; |w| = K,w\in L (y) \} .
$$
This completes the proof.\hfill $\Box$

\begin{lemme}
\label{sadiquelr}
Let $x$ be a primitive $S$-adic sequence on $A$ generated by $(\sigma_i : A_{i+1} \rightarrow A_i^* ; i\in \NN) \in S^{\pNN}$ and $a$ (with constant $s_0$). There exists a constant $K$ such that for all $b,c$ of $A_{s+1}$ and all integers $r,s$ with $s-r\geq s_0+1$ we have
$$
\frac{|\sigma_{r}\cdots \sigma_{s} (b)|}{|\sigma_{r}\cdots \sigma_{s}(c)|} \leq K.
$$
\end{lemme}

{\bf Proof.} For all non-negative integers $r$,  $s$ with  $r\leq s$ we set $S_{r,s} = \sigma_{r}\sigma_{r+1}\cdots \sigma_{s}$ and we denote by $M_{r,s}$ the matrix of this morphism. Let $r$ and $s$ be two non-negative integers such that $s-r\geq s_0 + 1$. The set $\{ S_{i,i+s_0} ; i\in \NN \}$ being finite we set
$$
K_1 = \max \{ |S_{i,i+s_0} (b) | ; i\in \NN,b\in A_{i+s_0+1} \} \ {\rm and} \ 
K_2 = \min \{ |S_{i,i+s_0} (b) | ; i\in \NN,b\in A_{i+s_0+1} \} .
$$
We remark that $K_2$ is not equal to $0$. Let $b$ and $c$ be two letters of $A_{s+1}$ and ${\bf 1} = (1,\cdots ,1)^T \in \NN^{|A|}$, we have 
$$
\frac{|S_{r,s}(b)|}{|S_{r,s}(c)|} = 
\frac{|S_{r,s-s_0-1}S_{s-s_0,s}(b)|}{|S_{r,s-s_0-1}S_{s-s_0,s}(c)|} \leq
\frac{||M_{r,s-s_0-1} (K_1 {\bf 1})||}{||M_{r,s-s_0-1} (K_2 {\bf 1})||} =
\frac{K_1}{K_2}
$$ 
where $|| (v_1 , \cdots , v_{|A|})^T || = \sum_{i=1}^{|A|} |v_i|$. This completes the proof.\hfill $\Box$

\medskip

{\bf Proof of Proposition \ref{sadiqlinrec}.} Suppose that $(X,T)$ is a primitive $S$-adic subshift generated by $(\sigma_i : A_{i+1} \rightarrow A_i^* ; i\in \NN)\in S^{\pNN}$ and $a$ (with constant $s_0$). Let $x = \lim_{n\rightarrow +\infty} \sigma_0\sigma_1\cdots \sigma_n (a)$. It follows from Lemma \ref{sadiquelb} that $x$ is uniformly recurrent and that the subshift  $(X,T)$ is minimal. We set $S_k = \sigma_{0}\cdots \sigma_{k}$ for all $k\in \NN$.  Let $u$ be a non-empty word of $L (X)$ such that $|u|\geq \max \{ |S_{s_0}(b)| ; b\in A_{s_0+1}\}$, and $v$ be a return word to $u$. We denote by $k_0$ the smallest positive integer $k$ such that $|u| < \min \{ |S_k (b)| ; b\in A_{k+1}\}$; we remark that $k_0\geq s_0+1$. We set $y = \lim_{n\rightarrow +\infty} \sigma_{k_0+1} \cdots \sigma_{n} (a)$. There exists a word of length 2, $bc$, of $L (y)$ such that $u$ has an occurrence in $S_{k_0} (bc)$. The sequence $y$ is uniformly recurrent (Lemma \ref{sadiquelb}) hence we can define $R$ to be the greatest difference between two successive occurrences in $y$ of a word of length 2 of $L (y)$. We have
$$
|v|
\leq R\max \{|S_{k_0} (d)| ; d\in A_{k_0+1} \} 
\leq RK \min \{|S_{k_0} (d)| ; d\in A_{k_0+1} \}
$$
$$
\leq RK \max \{|S_{k_0 -1} (d)| ; d\in A_{k_0} \} \min \{|\sigma_{k_0} (d)| ; d\in A_{k_0+1} \}
$$
$$
\leq RK^2 \min \{|S_{k_0 -1} (d)| ; d\in A_{k_0} \} \min \{|\sigma_{k_0} (d)| ; d\in A_{k_0+1} \}
$$
$$
\leq RK^2 \min \{|\sigma_{k_0} (d)| ; d\in A_{k_0+1} \} |u|,
$$
where $K$ is the constant of Lemma \ref{sadiquelr}. We set $M = RK^2 \min \{|\sigma_{i} (d)| ; i\in \NN , d\in A_{i+1} \} $. For all $u$ of $L (x)$ greater than $\max \{ |S_{s_0}(b)| ; b\in A_{s_0+1}\}$ and all $v$ in $\R_{u}$ we have $|v|\leq M|u|$. Hence $(X,T)$ is LR.

\vskip 0,3cm

We suppose now that $(X,T)$ is LR. The periodic case is trivial hence we suppose that $(X,T)$ is not periodic. There exists an integer $K$ such that 
$$
(\forall u\in L (X)) (\forall v\in \R_u ) ( \frac{1}{K} |u| \leq |v| \leq K |u| ),
$$
(this is Proposition \ref{proplr}). We set $\alpha = K^2(K+1)$. Let $x$ be an element of $X$. For each non-negative integer $n$ we set $\R_n = \R_{x,x_0x_1\cdots x_{\alpha^n-1}}$, $R_n = R_{x,x_0x_1\cdots x_{\alpha^n-1}}$ and $\Theta_n = \Theta_{x,x_0x_1\cdots x_{\alpha^n-1}}$. Let $n$ be a positive integer and $w$ be a return word to $x_0 x_1 \cdots x_{\alpha^n-1}$. The word $w$ is a concatenation of return words to $x_0 x_1 \cdots x_{\alpha^{n-1}-1}$. This induces a map $\lambda_n $ from $R_n $ to $R_{n-1}^{*} $ defined by  $\Theta_{n-1}\lambda_n = \Theta_{n}$. We set $ \lambda_0 = \Theta_0 $. For each letter $b$ of $R_n$ we have $|\Theta_{n-1}\lambda_n (b) | \leq K \alpha^n.$
Moreover each element of $\R_{n-1}$ is greater than $\alpha^{n-1}/K$ hence 
$$
|\lambda_n (b) | \leq \frac{K^2\alpha^n}{\alpha^{n-1}}\leq \alpha K^2.
$$
By Proposition \ref{proplr} we have $\card R_n \leq K(K+1)^2$, consequently the set $M = \{ \lambda_n ; n\in \NN^* \}$ is finite. The definition of $R_n$ implies that $\Theta_{n} (1)x_0 x_1\cdots x_{\alpha^n-1}$ is a prefix of $x$ for all $n\in \NN$ and 
$
\lambda_0 \lambda_1\cdots \lambda_n (1) = \Theta_{n} (1).
$
Lemma \ref{lemtech} implies that the length of $\Theta_{n} (1)$ tends to infinity with $n$ and
$$
x = \lim_{n\rightarrow + \infty} \lambda_0 \lambda_1 \cdots \lambda_n (11\cdots).
$$
Let $n\in \NN$. Each word of length $K\alpha^n$ has an occurrence in each word of length $(K+1)K\alpha^n$ (Proposition \ref{proplr}). Hence each element of $\R_n$ has an occurrence in each word of length $(K+1)K\alpha^n$. Let $w$ be an element of $\R_{n+1}$, we have $|w| \geq \alpha^{n+1}/K$ and $K(K+1)\alpha^n = \alpha^{n+1}/K $. Therefore each element of $\R_n$ has an occurrence in each element of $\R_{n+1}$. Consequently if $b$ belongs to $R_{n+1}$ then each letter of $R_n$ has an occurrence in $\lambda_n (b)$.

Consequently $(X,T)$ is a primitive $S$-adic subshift.\hfill $\Box$

\subsection{A necessary and sufficient condition for a sturmian subshift to be LR.}

We begin this subsection recalling some basic facts about sturmian subshifts (for more details see \cite{HM}). We prove that a sturmian subshift (generated by $\alpha \in \RR \setminus \QQ$) is LR if and only if the coefficients of the continued fraction of $\alpha$ are bounded.

Let $0<\alpha <1$ be an irrational number. We define the map $R_{\alpha}: [0,1[\rightarrow [0,1[$ by $R_{\alpha} (t) = t+\alpha$ (mod 1) and the map $I_{\alpha} : [0,1[\rightarrow \{ 0,1 \}$ by $I_{\alpha} (t) = 0$ if $t\in [0,1-\alpha[$ and $I_{\alpha} (t) = 1$ otherwise.  Let $\Omega_{\alpha} = \overline{ \{ (I_{\alpha} (R^n_{\alpha} (t)) ; n\in \ZZ) , t\in [0,1[ \} } \subset \{ 0,1 \}^{\pZZ}$. The subshift $(\Omega_{\alpha} , T)$ is called {\it sturmian subshift} (generated by $\alpha$) and its elements are called {\it sturmian sequences}.  There exists a factor map (see \cite{HM}) $\gamma : (\Omega_{\alpha} , T) \rightarrow ([0,1[, R_{\alpha})$ such that

\medskip

$\imath )$
$\card \gamma^{-1}( \{ \beta \})= 2$ if $\beta \in \{ n\alpha \; ({\rm mod} \; 1) ; n\in \ZZ \}$ and

$\imath \imath )$
$\card \gamma^{-1}( \{ \beta \})= 1$ otherwise. 

\medskip

Let $\beta \in [0,1[$ be an irrational number. It is well-known that $\Omega_{\alpha} = \Omega_{\beta}$ if and only if $\alpha = \beta$ and also that $(\Omega_{\alpha} , T)$ is a non-periodic minimal subshift (see \cite{HM}).

In what follows $\tau $ and $\sigma$ will be the morphisms from $\{ 0,1 \}$ to  $\{ 0,1 \}^*$ defined by 
\begin{center}
\begin{tabular}{lcl}
$\tau (0) = 0$  & and & $\sigma (0) = 01$ \\
$\tau (1) = 10$ &    & $\sigma (1) = 1 $.
\end{tabular}
\end{center}
An immediate consequence of a proposition (p. 206) of the  article \cite{AR} is the following.
\begin{prop}
Let $(X,T)$ be a subshift and $0<\alpha <1$ be an irrational number. Then $X = \Omega_{\alpha}$ if and only if $X = \Omega (x)$ where 
$$
x =
\lim_{k\rightarrow +\infty} \tau^{i_1} \sigma^{i_2} \tau^{i_3} \sigma^{i_4} \cdots \tau^{i_{2k-1}} \sigma^{i_{2k}} (00\cdots) .
$$
and $[0;i_1 + 1,i_2,\cdots ]$ is the continued fraction of $\alpha $.
\end{prop}

Therefore each sturmian sequence generates a $S$-adic subshift. From this proposition we can characterize those sturmian subshifts which are LR.
\begin{prop}
A sturmian subshift $(\Omega_{\alpha}, T)$ is LR if and only if the coefficients of the continued fraction of $\alpha$ are bounded.
\end{prop}
{\bf Proof.} First we remark that $\tau^n (1) = 10^n$ and $\sigma^n (0) = 01^n$ for all $n\in \NN$. Let $[0;i_1+1,i_2,\cdots]$ be the continued fraction of $\alpha$. If the sequence $(i_n ; n\in \NN)$ is not bounded then for all $n\in \NN$ there would exist a word $u$ such that $u^n\in L(\Omega_{\alpha})$. Thus $(\Omega_{\alpha} , T)$ could not be LR (Proposition \ref{proplr}).

If the sequence $(i_n ; n\in \NN )$ is bounded, we can check that $(\Omega_{\alpha} , T)$ is a primitive $S$-adic subshift. Proposition \ref{sadiqlinrec} achieves the proof.
\hfill $\Box$

\medskip

Consequently it follows from Theorem \ref{maintheo} that if the coefficients of the continued fraction of $\alpha$ are bounded then $(\Omega_{\alpha} , T)$ has a finite number of subshift factors up to isomorphism. This is not surprising because sturmian subshifts do not have proper Cantor factor (Corollary \ref{facteursturmfini}) although they have infinitely many factors which are not Cantor systems. 

\subsection{Sturmian subshifts are Cantor prime.}

In this subsection we prove that sturmian subshifts are {\it Cantor prime}, i.e. they do not have proper Cantor factor. It follows from the following proposition.
\begin{prop}
Let $0<\alpha <1$ be an irrational number. Let $(X,S)$ be a non-periodic factor of $(\Omega_{\alpha} , T)$ which is not isomorphic to $(\Omega_{\alpha} , T)$. Then $(X,S)$ is a factor of $([0,1[ , R_{\alpha})$.
\end{prop}
{\bf Proof.} Let $\O = \{ n\alpha \; ({\rm mod} \; 1) ; n\in \ZZ \}$. We set $\gamma^{-1}( \{ \beta \}) = \{ x(\beta), y(\beta) \}$ for all $\beta \in \O$ and $\gamma^{-1}( \{ \beta \}) = \{ x(\beta) \}$ for all $\beta \in [0,1[ \setminus \O$ .

Let $\phi : (\Omega_{\alpha} , T) \rightarrow (X,S)$ be a factor map.  It suffices to prove that $\phi (x(\beta)) = \phi (y(\beta))$ for all $\beta \in \O$; i.e. that there exists $\beta^{'}\in \O$ such that $\phi (x(\beta^{'})) = \phi (y(\beta^{'}))$. Because if we set $\rho (\beta ) = \phi (x(\beta))$, for all $\beta \in [0,1[$ it would define a factor map $\rho$ from $([0,1[, R_{\alpha})$ onto $(X,S)$.

The map $\phi$ is not an isomorphism hence there exist distinct elements $x^1,x^2 \in X$ such that $\phi (x^1) = \phi (x^2)$. We consider two cases.

First case: $ \gamma (x^2) -  \gamma (x^1) \in \O$; there exists $k\in \ZZ$ such that $\gamma (T^k x^1) = \gamma (x^2)$.

Let $(n_i ; i\in \NN )$ be a sequence of integers such that $(T^{n_i} x^1 ; i\in \NN )$ converges to $y\in \Omega_{\alpha} \setminus \gamma^{-1} (\O ) $. We can suppose that $(T^{n_i} x^2 ; i\in \NN )$ converges to some $z$. Then $\gamma (T^k y) = \gamma (z)$ and $T^k y = z$ because $y\not \in \gamma^{-1} (\O )$. On the other hand we obtain $\phi (y) = \phi (z) = \phi (T^k y)$. But $(X,S)$ is not periodic thus $k=0$; i.e $\gamma (x^1) = \gamma (x^2)$ and consequently $x^1\in \O$.

Second case: $ \gamma (x^2) -  \gamma (x^1) \not \in \O$. 

Let $(m_i ; i\in \NN )$ and $(n_i ; i\in \NN )$ be sequences of integers such that $(T^{m_i} x^1 ; i\in \NN )$ and $(T^{n_i} x^1 ; i\in \NN )$ converge respectively to $x(0)$ and $y(0)$. We can suppose that $(T^{m_i} x^2 ; i\in \NN )$ and $(T^{n_i} x^2 ; i\in \NN )$ converge respectively to $z^1$ and $z^2$. We remark that $\gamma ( z^1 ) = \gamma (x^2) -  \gamma (x^1) = \gamma ( z^2 )$. As $ \gamma (x^2) -  \gamma (x^1) \not \in \O$  we deduce that $z^1 = z^2$. It follows that $\phi ( x(0)) = \phi (z^1) = \phi (z^2) = \phi (y(0))$ which ends the proof.
\hfill $\Box$

\medskip

It is well-known and easy to prove that if $([0,1[ , R_{\alpha})$ has a $p$-periodic factor then $p=1$. Moreover $([0,1[ , R_{\alpha})$ cannot have a non-periodic Cantor factor for topological reasons. Hence we obtain

\begin{coro}
\label{facteursturmfini}
Let $\alpha \in [0,1[$ be an irrational number. Then $(\Omega_{\alpha} , T)$ is Cantor prime.
\end{coro}
We recall that sturmian subshifts $(\Omega_{\alpha} , T)$ are not prime because all rotations $([0,1[ , R_{n\alpha})$ are factors of $(\Omega_{\alpha} , T)$.

\section{LR subshifts are uniquely ergodic.}

Let $(X,S)$ be a dynamical system. An {\it invariant measure} for $(X,S)$ is a probability measure $\mu$, on the $\sigma$-algebra $\B (X)$ of Borel sets, with $\mu (S^{-1} B) = \mu (B)$ for all $B\in \B (X) $; the measure is {\it ergodic} if every $S$-invariant Borel set has measure 0 or 1. The set of invariant measures for $(X,S)$ is denoted by $\M (X,S)$. The system $(X,S)$ is {\it uniquely ergodic} if $\card (\M (X,S)) = 1$.
\begin{prop}
\label{encadmes}
Let $(X,T)$ be a LR subshift with constant $K$. Then for all $\mu \in \M (X,T)$ and all $u\in L(X)$ we have
$$
1/K \leq |u| \mu ([u]) \leq K  .
$$
\end{prop}
{\bf Proof.} It suffices to prove the result for ergodic measures. Let $\mu \in \M (X,T)$ be an ergodic measure. Let $u\in L (X)$ and $\xi$ be the characteristic function of the cylinder $[u]$. We remark that for all $x\in X$ we have
$$
\frac{1}{n} \sum_{i=0}^{n-1} \xi (T^i x) = \frac{L_u(x_{[0,n + |u|-2]})}{n} \leq \frac{1}{n} \frac{n + |u|}{(1/K)|u|}  .
$$
We can apply the Ergodic Theorem (see \cite{Wa}) to obtain that for $\mu$-almost every $x\in X$,
$$
\lim_{n\rightarrow +\infty} \frac{1}{n} \sum_{i=0}^{n-1} \xi (T^i x) = \mu ([u]) 
$$
Hence $|u| \mu ([u]) \leq K$. In the same way we obtain $(1/K) \leq |u| \mu ([u])$.\hfill $\Box$ 

\medskip

In \cite{Bo} Boshernitzan proved the following result.
\begin{theo}
\label{bosh}
Let $(X,T)$ be a minimal subshift which is not uniquely ergodic. Let $\mu$ be a invariant measure for $(X,T)$, then 
$$
\lim_{n\rightarrow +\infty} n \epsilon (n) = 0 \ where \ \epsilon (n) = \min \{ \mu ([u]) ; u\in L_n (X)  \}.
$$
\end{theo}
It follows directly from Proposition \ref{encadmes} and Theorem \ref{bosh} that :
\begin{theo}
Linearly recurrent subshifts are uniquely ergodic.
\end{theo}

\label{sectiontrois}

\section{Factors of LR subshifts.}
\label{sectionquatre}

In this section we prove Theorem \ref{maintheo}.
\subsection{Non-periodic LR subshifts are coalescent.}
\label{coalsection}

A dynamical system $(X,S)$ is {\it coalescent} if each factor map $\gamma : (X,S) \rightarrow (X,S)$ is an isomorphism. Two dynamical systems $(X,S)$ and $(Y,S')$ are {\it weakly isomorphic} if each is a factor of the other. 

In \cite{DHS} the following theorem is proved.
\begin{theo} 
\label{chainefini}
Let $(Y,T)$ be a non-periodic LR subshift. There exists a constant $D$ such that if 
$$ 
(Y,T)\; \fleche{ \gamma_{D-1}} \;
(Y_{D-1},T)\;  \fleche{ \gamma_{D-2}}\; 
\cdots\; \fleche{\gamma_1}(Y_1,T)\;
\fleche{\gamma_0} \; (Y_0,T)\  ,
$$
is a chain of non-periodic subshift factor maps then there exists $0\leq n\leq D-1$ such that $\gamma_{n}$ is an isomorphism.
\end{theo}
It follows that 
\begin{coro}
Let $(X,T)$ and $(Y,T)$ be two {\it weakly isomorphic} LR subshifts. Then $(X,T)$ and $(Y,T)$ are isomorphic.
\end{coro}
{\bf Proof.}
The case of periodic subshifts is trivial, hence we suppose that $(X,T)$ and $(Y,T)$ are non-periodic.

Let $\phi : (X,T)\rightarrow (Y,T)$ and $\psi : (Y,T)\rightarrow (X,T)$  be two factor maps and consider the following chain of factors
$$ 
(X,T)\; \fleche{ \phi } \;
(Y,T)\;  \fleche{ \psi}\; 
(X,T)\; \fleche{ \phi } \;
\cdots\; \fleche{\psi }(X,T)\;
\fleche{\phi} \; (Y,T)\  ,
$$
It suffices to apply twice Theorem \ref{chainefini} to prove that $\phi$ and $\psi$ are isomorphism which achieves the proof.
\hfill $\Box$

\medskip

It follows directly from the proof of the previous corollary that
\begin{coro}
Non-periodic LR subhifts are coalescent.
\end{coro}

\subsection{Preimages of factor maps of LR subshifts.}

Let $(X,S)$ and $(Y,T)$ be two dynamical systems. The factor map $\phi : (X,S) \rightarrow (Y,T)$ is {\it finite-to-one} (with constant $K$) if for all $y\in Y$ we have $\card (\phi^{-1} (\{ y \} )) \leq K$.

Let $\phi$ be a factor map from the subshift
$(X,T)$ on the alphabet $A$ onto the subshift $(Y,T)$ on the
alphabet $B$. The theorem of Curtis-Hedlund-Lyndon (Theorem 6.2.9 in \cite{LM}) asserts
that $\phi$ is a {\it sliding block code}. That is to
say there exists a {\it $r$-block map} $f: A^{2r+1}
\rightarrow B$ such that $(\phi (x))_i = f(x_{[i-r,i+r]})$ for
all $i\in \ZZ$ and $x\in X$. We shall say that $f$ is a {\it block map
associated to $\phi$} and that $f$ {\it defines} $\phi $. If $u= u_0 u_1 \cdots u_{n-1}$ is a word of length $n\geq 2r+1$ we define $f(u)$ by $(f(u))_i = f(u_{[i,i+2r]})$, 
$i\in \{ 0,1,\cdots , n-2r-1 \}$. Let $C$ denote the alphabet $A^{2r+1}$ and $Z=\{
((x_{[-r+i,r+i]}) ; i\in \ZZ)\in C^{\pZZ }; (x_n ; n\in  
\ZZ ) \in X \}$. It is easy to check that the subshift $(Z,T)$
is isomorphic to $(X,T)$ and that $f$ induces a 0-block map from $C$
onto $B$ which defines a factor map from $(Z,T)$ onto $(Y,T)$.

\begin{lemme}
\label{mainlemme}
Let $(X,T)$ be a non-periodic LR subshift (with constant $K$) and $(Y,T)$ be a non-periodic subshift factor of $(X,T)$. Then
 $(Y,T)$ is LR. Moreover there exists $n_1\in \NN$ such
that: 
For all $u\in L(Y) $, with $|u|\geq n_1$, we have
\begin{enumerate}
\item
\label{point1}
$|u|/2K \leq |w|\leq 2K|u|$ for all $w\in  \R_u$ ; 
\item
\label{point2}
$\card (\R_u ) \leq 2K(2K+1)^2$.
\end{enumerate}
\end{lemme}
{\bf Proof.} Point \ref{point1} is proved in \cite{DHS} and gives the constant $n_1$. Let $u\in L(Y)$ with $|u|\geq n_1$ and $v \in L(Y)$ be a word of length $(2K+1)^2|u|$. Each word of length $(2K+1)|u|$ occurs in $v$,
hence each return word to $u$ occurs in $v$. It follows from 
the previous assertion that in $v$ occurs at the most 
$2K(2K+1)^2|u|/|u|= 2K(2K+1)^2$ return words to $u$.
\hfill $\Box$

\begin{theo}
Let $(X,T)$ be a non-periodic LR subshift (with constant $K$). If $\phi : (X,T) \rightarrow (Y,T)$ is a factor map such that $(Y,T)$ is a non-periodic subshift then $\phi$ is finite-to-one with constant $4K(K+1)$.
\end{theo}

{\bf Proof.}  Let $\phi : (X,T) \rightarrow (Y,T)$ be a factor map where $(Y,T)$ is a non-periodic subshift on the alphabet $B$. Let $f : A^{2r+1} \rightarrow B$ be a block map defining $\phi$. It suffices to prove that there exists an integer $n_0$ such that for all $u\in f ( L(X))$, with $|u| \geq n_0$, we have $\card ( f^{-1} (\{ u \} ))\leq 4K(K+1)$.

Let $n_1$ be the integer given by Lemma \ref{mainlemme}. We set $n_0 = \max (2r +1, n_1)$. Let $u\in f(L(X))$ such that $|u| \geq  n_0$. The difference between two distinct occurrences of elements of $ f^{-1} (\{ u \} )$ is greater than $|u|/2K $. Moreover $ f^{-1} (\{ u \} ) \subset L_{|u|+2r} (X) $ and each word of length $(K+1)(|u|+2r)$ has an occurrence of each word of $L_{|u|+2r} (X)$. Therefore
$$
\card ( f^{-1} (\{ u \} )) \leq \frac{(K+1)(|u|+2r)}{|u|/2K} \leq 4K (K+1).
$$ 
This completes the proof.
\hfill $\Box$

\subsection{Proof of Theorem \protect{\ref{maintheo}}.}

We first begin with a lemma.
\begin{lemme}
\label{lemmeffini}
Let $(X,S)$ be a minimal dynamical system and $\phi_1 : (X,S) \rightarrow (X_1,S_1)$, $\phi_2 : (X,S) \rightarrow (X_2,S_2)$ be two factor maps where $(X_1,S_1)$ and $(X_2,S_2)$ are non-periodic and where $\phi_1$ is finite-to-one. If there exist $x,y\in X$ and $r\in \ZZ$ such that $\phi_1 (x) = \phi_1 (y)$ and $\phi_2 (x) = S_2^r \phi_2 (y)$, then $r=0$.
\end{lemme}
{\bf Proof.} Let $y_1 = y$. There exists $(n_i ; i\in \NN )$ such that $\lim_{i\rightarrow +\infty } S^{n_i} x=y$. By compactness we can suppose that $\lim_{i\rightarrow +\infty } S^{n_i} y$ exists, we call it $y_2$. We have  $\phi_1 (y_1) = \phi_1 (y_2)$ and $\phi_2 (y_1) = S_2^r \phi_2 (y_2)$. By compactness we can suppose that $\lim_{i\rightarrow +\infty } S^{n_i} y_2$ exists, we call it $y_3$. Thus we have $\phi_1 (y_2) = \phi_1 (y_3)$ and $\phi_2 (y_2) = S_2^r \phi_2 (y_3)$. Hence we can construct a sequence $(y_n ; n\in \NN )$ of elements of $X$ such that $\phi_1 (y_i) = \phi_1 (y_{i+1})$ and $\phi_2 (y_i) = S_2^r \phi_2 (y_{i+1})$ for all $i\geq 1$. But $\phi_1$ is finite to one, therefore there exist $i<j$ such that $y_i = y_j$. Then we have 
$$
\phi_2 (y_i) = S_2^r \phi_2 (y_{i+1}) =  S_2^{2r} \phi_2 (y_{i+2}) = \cdots = S_2^{r(j-i)} \phi_2 (y_{j}) = S_2^{r(j-i)} \phi_2 (y_{i}).
$$
Consequently $r=0$ because $(X_2,S_2)$ is non-periodic.\hfill $\Box$

\medskip

We can now prove Theorem \ref{maintheo}.

Let $(X,T)$ be a LR subshift, on the alphabet $A$, with constant $K$. Suppose that the set of subshift factors of $(X,T)$ is infinite (up to isomorphism) and let $\F= \{ \phi_i : (X,T) \rightarrow (X_i,T) ; 1\leq i \leq N \}$, with $N > ((2K(2K+1)^2)^{4K^2})^{K(K+1)^2} $, be a set of subshift factor maps such that for all $1\leq i<j \leq N$ the subshifts $(X_i,T)$ and $(X_j , T)$ are not isomorphic. Let $A_i$, $ 1\leq i \leq N $, be the alphabet of $(X_i ,T)$. By the Curtis-Hedlund-Lyndon Theorem there exists a positive integer $r$ such that for all $i\in \{1,\cdots ,N \}$ there is a map $f_i : A^{2r+1} \rightarrow A_i$ satisfying 
$$
(\phi_i (x))_j = f_i (x_{[ j-r,j+r]}) \ \ \forall x\in X , \ \ \forall j\in \ZZ .
$$
Considering $B = A^{2r+1}$ as an alphabet, the system $(Y,T)$, where  $Y = \{ ((x_{n-r},x_{n+r});n\in \ZZ) ; (x_n;n\in \ZZ ) \in X \}$, is a subshift on $B$ which is isomorphic to $(X,T)$. For all $i\in \{1,\cdots , N \}$ let $g_i : B\rightarrow A_i$ be the map defined for all $(u)\in B$  by $g_i ((u)) = f_i(u) $ and $\psi_i : (Y,T) \rightarrow (X_i,T)$ be the factor map defined by $g_i$. 

Due to Lemma \ref{mainlemme} there exists $n_0\in \NN$ such that for all $u\in L(Y) $ with $|u| \geq n_0$ we have 
\begin{itemize}
\item
$|u|/2K \leq |w|\leq 2K|u| $ for all $w\in  \R_u$ ;
\item
$|g_i (u)|/2K \leq |w|\leq 2K|g_i (u)| $ and $\card (\R_{g_i (u)}) \leq 2K (2K+1)^2$ for all $i \in \{ 1 , \cdots , N \}$ and all $w\in  \R_{g_i(u)}$.
\end{itemize}
Let $u$ be a word of $L(Y)$ such that $|u|\geq n_0$ and $i$ be an element of $ \{ 1 , \cdots , N \}$. We set $u_i = g_i(u)$. For each return word $w$ to $u$ the word $g_i ( w)$ is a concatenation of return words to $u_i$. This induces a map $\lambda_i$ from $R_u$ to $R_{u_i}^*$ defined by $g_i \Theta_u = \Theta_{u_i} \lambda_i$. This map is such that
$$
|\lambda_i (b)| \leq \frac{2K|u|}{|u_i|/2K} = 4K^2.
$$
Moreover from  Proposition \ref{proplr} we know that $\card (\R_v) \leq K (K+1)^2 $ for all $v\in L(X)$. Hence we obtain $\card \{ \lambda_i ; 1\leq i \leq N \} \leq ((2K(2K+1)^2)^{4K^2})^{K(K+1)^2} $. But $N$ is strictly greater than $ ((2K(2K+1)^2)^{4K^2})^{K(K+1)^2} $, thus there exist $i$ and $j$, with $i\not = j$, such that $\lambda_i = \lambda_j$. We set $\lambda = \lambda_i$.

\medskip

We claim that for all $x,y$ in $Y$, if $\psi_i (x) = \psi_i (y)$ then $\psi_j (x) = \psi_j (y)$. First we remark that it suffices to prove this claim for all $x\in [u]$ and all $y\in Y$.

Let $x\in [u]$ and $y \in Y$ such that $\psi_i (x) = \psi_i (y)$. Let $r\in \NN$ be such that $T^r y \in [u]$. There exists a unique sequence $\xb\ \in R_u^{\pZZ}$ (resp. $\yb \in R_u^{\pZZ}$) such that $\Theta_u (\xb ) = x$ (resp. $ \Theta_u (\yb ) = T^r y$). We set $\Theta_i = \Theta_{u_i}$ and $\Theta_j = \Theta_{u_j}$. We have
$$
\psi_i (x) = \psi_i \Theta_u (\xb) = \Theta_i \lambda (\xb) \ \hbox{\rm and} \ \psi_i (y) = T^{-r} \psi_i \Theta_u (\yb) = T^{-r} \Theta_i \lambda (\yb).
$$
But $\psi_i (x) = \psi_i (y)$ consequently $T^{-r} \psi_i \Theta_u (\yb) \in [u_i]$ where $r = | \Theta_i (m) |$, for some prefix $m$ of $\lambda (\xb )^+$, and $T^{|m|} \lambda (\xb ) = \lambda (\yb )$. On the other hand we have:
$$
\psi_j (x) = \psi_j \Theta_u (\xb) = \Theta_j \lambda (\xb) \ \hbox{\rm and} \ \psi_j (y) = T^{-r} \Theta_j \lambda (\yb) = T^{-r} \Theta_j T^{|m|} \lambda (\xb ) = T^{|\Theta_j (m) | - r } \Theta_j \lambda (\xb ).
$$
Hence we have $\psi_j (y) = T^{ |\Theta_j (m) | - |\Theta_i (m) | } \psi_j (x) $ and Lemma \ref{lemmeffini} implies $\psi_j (y) = \psi_j (x)$. The claim is proved.

Consequently we can define the map $\psi : X_i \rightarrow X_j$, for all $x\in X_i$ by $\psi (x) = \psi_j(x')$ where $x'$ belongs to $\psi_i ^{-1} ( \{ x \} )$. There is no difficulty to prove that $\psi $ is an isomorphism from $X_i$ onto $X_j$. This is in contradiction with the definition of $\F$.\hfill $\Box$

\begin{coro}
\label{maincoro}
Substitution subshifts have a finite number of subshift factors up to isomorphism.
\end{coro}
{\bf Proof.} In \cite{DHS} it is proved that substitution subshifts are LR.\hfill $\Box$

\section{Cantor factors of substitution subshifts.}
\label{sectioncinq}

The goal of this section is to give a necessary and sufficient condition for a substitution system to have a finite number of non-periodic Cantor factors.

Let $(X,T)$ be a substitution subshift, $\F (X,T)$ (resp. $\F^{*} (X,T)$) will denote the set of its (resp. non-periodic) Cantor factors. In the sequel when we will say that $\F (X,T)$ or $\F^* (X,T)$ is finite it will mean that it is finite up to isomorphism. In \cite{DHS} the following proposition is proved.
\begin{prop}
\label{subouodo}
The Cantor factors of substitution subshifts are substitution subshifts or odometers with constant base. 
\end{prop}

We need the following notion. For a dynamical system $(Z,S)$ the {\it periodic spectrum of $(Z,S)$}, $\P (Z,S)$, is the set of positive integers $p$ for which there are disjoint clopen sets $Y,  \ SY, \ \cdots , S^{p-1} Y $ whose union is $Z$. We remark that the odometer with base $(p,p,\cdots)$ is a factor of $(Z,S)$ if and only if for all $n\in \NN$ $p^n$ belongs to $\P (Z,S)$.

\begin{lemme}
\label{lemspec}
Let $(X,T)$ be a non-periodic substitution subshift. Then

$\imath )$
The set of prime numbers belonging to $\P (X,T)$ is finite.

$\imath\imath )$
$\F (X,T)$ is finite if and only if $\P (X,T)$ is finite.

$\imath\imath\imath )$ $\F^{*} (X,T)$ is finite if and only if it does not exist two distinct prime numbers $p,q\in \P (X,T)$ such that $p^n,q^n\in \P (X,T)$ for all $n\in \NN$.
\end{lemme}
{\bf Proof.} We will only prove $\imath)$; $\imath\imath)$ and $\imath\imath\imath)$ are left to the reader. 

Let $\sigma : A \rightarrow A^*$ be a substitution generating $(X,T)$, $M$ be its matrix and $a\in A$ be such that the first letter of $\sigma (a)$ is $a$. Let $p$ be a prime number and $Y$ be a clopen set such that $Y, TY,\cdots , T^{p-1} Y$ is a clopen partition of $X$. Let $u$ be a return word to $a$. There exists $n_0$ such that for all $n\geq n_0$ $p$ divides $|\sigma^n (u)|$. Let $n\geq n_0$ and $P(X) = \sum_{i=0}^{|A|} c_i X^i$ be the characteristic polynomial of $M^n$. We know that $c_0 = \det M^n$ and by the Theorem of Cayley-Hamilton that $P ( M^n) u = 0$. Therefore $p$ divides $|u|\det M$ and that $\P (X,T)$ is finite.
\hfill $\Box$

\medskip

{\bf The constant length case.} For the constant length case the necessary and sufficient condition is easy to obtain on account of a result due to Dekking \cite{De} which gives a complete characterization of the maximal equicontinuous factors of subshifts generated by substitutions of constant length.

\begin{theo}
{\rm \cite{De}}
\label{dek}
Let $(X,T)$ be a subshift generated by a substitution of constant length $l$. The maximal equicontinuous factor of $(X,T)$ is an odometer with base $(h,l,l,l, \cdots )$ where $h = \max \{ n\geq 1 ; (n,l) = 1 , n \ \hbox{\rm divides } \ \gcd \{ i\geq 0 ; x_i = x_0 \} \}$, where $x\in X$. 
\end{theo}
Odometers with base $(p_n \geq 2 ; n\in \NN )$ have an infinite number of non-isomorphic periodic factors, hence if $(X,T)$ is non-periodic and generated by a substitution of constant length then $\F (X,T)$ is infinite. But it can have a finite number of non-periodic factors, this is the next result.
\begin{coro}
Let $(X,T)$ be a non-periodic subshift generated by a substitution of constant length $l$. 

$\imath )$ $\F (X,T)$ is infinite.

$\imath\imath )$ $\F^{*} (X,T)$ is finite if and only if $l$ is a prime number. 
\end{coro}
{\bf Proof.} It follows from Corollary \ref{maincoro}, Proposition \ref{subouodo}, Lemma \ref{lemspec} and Theorem \ref{dek} and the previous remarks. \hfill $\Box$

\medskip

{\bf General case.} The following result is proved in \cite{DHS}.

\begin{lemme}
Let $(X,T)$ be a non-periodic substitution subshift generated by $\sigma : A\rightarrow A^*$. For every $n>0$,
$$
\P_n = \{ T^k (\sigma^n ([b])) ; b\in A , 0\leq k\leq |\sigma^n (b)| - 1 \}
$$
is a clopen partition of $X$.
\end{lemme}

A substitution on the alphabet $A$ is {\it proper} if there exists $a \in A$ such that for all $c\in A$ the word $\sigma (c)$ begins with the letter $a$. In \cite{DHS} the definition is slightly different.

\begin{lemme}
\label{reco}
Let $(X,T)$ be a non-periodic subshift generated by a proper substitution, $\sigma : A \rightarrow A^*$, and $p$ be an integer.  The following are equivalent.

$\imath )$ For all $n\in \NN$ there exists an integer $m$ such that $p^n$ divides $\gcd \{ |\sigma^m (c)| ; c\in A \}$.

$\imath \imath )$ $p^n$ is in $\P(X,T)$ for all $n\in \NN$.
\end{lemme}
{\bf Proof.} Let $a$ be a letter of $A$ such that for all $c\in A$ the word $\sigma (c)$ begins with the letter $a$. Let $n,m\in \NN$ be such that $p^n$ divides $\gcd \{ |\sigma^m (c)| ; c\in A \}$.

We know that 
$$
\P_m = \{ T^k (\sigma^m ([c])) ; c\in A , 0\leq k\leq |\sigma^m (c)| -1 \}
$$
is a clopen partition of $X$. Consequently $X = \cup_{i=0}^{p^n-1} T^i Y$ where 
$$
Y = \bigcup_{c\in A} \bigcup_{j=0}^{(|\sigma^m (c)|/p^n)-1} T^{jp^n} \sigma^m ([c]),
$$
i.e $p^n$ is in $\P(X,T)$.

Suppose that $p^n$ is in the periodic spectrum of $(X,T)$ for all $n\in \NN$. There exists a clopen set $Y$ such that $Y, TY, \cdots,T^{p^n-1}Y$ is a clopen partition of $X$. Let $uv\in L (X)$ such that $[u.v]\subset Y$. The length of the return words to $uv$ are in $p^n \ZZ$. Let $m$ be such that $uv$ has an occurrence in $\sigma^{m-1} (a)$. For all $c\in A$ the word $\sigma^m (c)$ is a concatenation of return words to $\sigma^{m-1} (a)$, consequently $|\sigma^m (c)| \in p^n \ZZ$. 
\hfill $\Box$

\medskip

The following lemma determines the set of prime numbers satisfying Point $\imath\imath )$ of Lemma \ref{reco}.
\begin{lemme}
\label{lemmemat}
Let $M$ be a $d\times d$ integral matrix and $e\in \ZZ^d$ be the vector which coefficients are equal to 1. Let $p$ be a prime number, the following are equivalent:
\begin{enumerate}
\item
\label{cond1}
$\forall n \in \NN, \ \exists k\in \NN, \ M^ke \in p^n \ZZ^d  $,
\item
\label{cond2}
$p$ divides $\gcd (a_0 , \cdots , a_r )$ where $r = \max \{ i\in \NN ; \{ e, Me,\cdots, M^i e \} \ \ \rm{is \ free} \}$ and $Q(X) = \sum_{i=0}^{r+1} a_i X^i \in \ZZ [X]$ is the characteristic polynomial of the restriction of $M$ to the vector subspace spanned by $e,Me,\cdots,M^r e$.

\end{enumerate}
\end{lemme}
{\bf Proof.} Let $P(X) = \sum_{i=0}^{d} \alpha_i X^i$ be the characteristic polynomial of $M$. We first prove that \ref{cond1} implies \ref{cond2}.

First case: $r = d-1$.

We have $P(X) = Q(X)$. Let $f$ be a vector of the canonical base of $\ZZ^d$. There exist $p_0/q_0, \cdots ,p_{d-1}/q_{d-1} \in \QQ$ such that $f = \sum_{i=0}^{d-1} (p_i/q_i) M^i e$. Let $n = \max \{ i \in \NN ; \exists j\in \{ 0,1, \cdots, d-1\} , q_j \in p^i \ZZ \  \}$ and $k$ be such that $M^k e \in p^{n+1} \ZZ^d $ (this implies that $M^{k+i} e \in p^{n+1} \ZZ^d $ for all $i\in \NN$). We have $M^k f = \sum_{i=0}^{d-1} (p_i/q_i) M^{k+i} e$, consequently $M^k f$ belongs to $p\ZZ^d$. It follows that for all $i\in \NN$ the coefficients of $M^{k+i}$ belongs to $p\ZZ$.

Let $l\in \NN$ be such that $p^l\geq k$ and $R(X) = \sum_{i=0}^d b_i X^i$ be the characteristic polynomial of $M^{p^l}$. We notice that $p$ divides $\gcd (b_0, \cdots , b_{d-1})$. The number $p$ is prime therefore we obtain (in $(\ZZ/p\ZZ) [X]$)
$$
\sum_{i=0}^d a_i^{p^l} X^{ip^l} = (\sum_{i=0}^{d} a_i X^i)^{p^l} = \det (M-X\id)^{p^l} = \det (M^{p^l} - X^{p^l}\id ) = R( X^{p^l}) = X^{dp^l}.
$$
Thus $p$ divides $\gcd (a_0^{p^l}, \cdots , a_{d-1}^{p^l})$ and clearly $p$ divides $\gcd (a_0, \cdots , a_{d-1})$.

\medskip

Second case: $r\leq d-2$.

Let $N$ be the $(r+1)\times(r+1)$ matrix which is the restriction of $M$ to the vector subspace spanned by $e,Me,\cdots , M^re$. It has the following form
\[
N= \left[
\begin{array}{cccccc}
0       & 0       & \cdots & 0      & c_0     \cr 
1       & 0       & \ddots & \vdots & c_1     \cr
0       & 1       & \ddots & 0      & c_2     \cr
\vdots  & \ddots  & \ddots & 0      & \vdots  \cr
0       & \cdots  & 0      & 1      & c_r     \cr
\end{array}
\right] ,
\]
where $c_i \in \QQ$ for all $0\leq i\leq r$, and consequently $Q(X) = \sum^{r+1}_{i=0} c_i X^i$ where $c_{r+1} = (-1)^{r+1}$.

There exists $R(X) = \sum_{i=0}^{d-r-1} (u_i/v_i) X^i \in \QQ [X]$ such that $P(X) = Q(X) R(X)$. We set $c_i = r_i/s_i$, $0\leq i \leq r+1$, $s = \scm ( s_0 ,\cdots , s_{r+1})$ and $v = \scm ( v_0, \cdots , v_{d-r-1} )$. We have $Q(X) = Q^{'} (X)/s$ and $R(X) = R^{'} (X)/v$ where $Q^{'} (X)$ and $R^{'} (X)$ belong to $\ZZ [X]$. Given $F(X)\in \ZZ[X]$ we denote by $c(F)$ the great common divisor of the coefficients of $F$. We easily see that $c(Q^{'})$ (resp. $c(R^{'})$) divides $s$ (resp. $v$); let $s^{'}$ and $v^{'}$ be respectively the integers $s/c(Q^{'})$ and $v^{'} =  v/c(R^{'})$. Applying Gauss Lemma we obtain $sv c (P) = c(Q^{'})c(R^{'})$ and then $s^{'}v^{'}=1$ because $c(P) = 1$. Thus $Q(X) \in \ZZ[X]$ and a fortiori $c_i \in \ZZ, 0\leq i\leq r$.

The matrix $N$ fulfills the condition of the first case, consequently $p$ divides $\gcd (c_0, \cdots , c_r)$.

\medskip

To prove that \ref{cond2} implies \ref{cond1} it suffices to notice that for all $i\in \NN$  we have
$$
-a_{r+1} M^{r+1+i}e = a_r M^{r+i} e + a_{r-1} M^{r-1+i} e + \cdots + a_0 M^i e,
$$
where $a_{r+1} = 1$ or $-1$.\hfill $\Box$

\begin{prop}
Let $(X,T)$ be a non-periodic substitution subshift generated by a proper substitution which transposed matrix is $M$. Let $Q(X) = \sum_{i=0}^{r+1} a_i X^i$ be the polynomial given by Lemma \ref{lemmemat}. Then 

$\imath )$
The set $\F (X,T)$ is finite if and only if $\gcd (a_0,\cdots , a_r) = 1$. 

$\imath\imath )$ The set of $\F^* (X,T)$ is finite if and only if $\gcd (a_0,\cdots , a_r)$ is a prime number.
\end{prop}
{\bf Proof.} It follows from Corollary \ref{maincoro}, Proposition \ref{subouodo} and Lemmas \ref{lemspec}, \ref{reco} and \ref{lemmemat}.
\hfill $\Box$

\medskip

It remains the case of non proper substitutions.
\begin{prop}
\label{subprop}
Let $(X,T)$ be a substitution subshift generated by $\sigma :A\rightarrow A^*$. Then there exists a proper substitution $\zeta :B\rightarrow B^* $ such that the subshift it generates is isomorphic to $(X,T)$. (Moreover the substitution $\zeta$ can be computed explicitely.)
\end{prop}
{\bf Proof.}
This is proved in \cite{DHS} (Section 6). We just recall the construction of $\zeta$. 

Let $a\in A$ be such that $\sigma (a)$ begins with $a$. Let $x$ be the fixed point of $\sigma$ such that $x_0 = a$. We recall that $\R_a$ is the set of return words to $a$ and that $\Theta = \Theta_a : R_a \rightarrow \R_a$ is the map defined in Section \ref{sectionun}. Let $\tau : R_a \rightarrow R_a^*$ be the only morphism such that
$$
\sigma \Theta = \Theta \tau.
$$
The substitution $\tau$ has the following properties:

\medskip

$\imath )$ $\tau (b)$ begins with the letter 1 for all $b\in R_a$, i.e $\tau $ is proper,

$\imath \imath )$ $\tau$ has a unique fixed point $y$ and it satisfies $\Theta (y) = x$.

\medskip

For the details we refer the reader to \cite{Du1,Du2}.

Substituting a power of $\tau$ for $\tau$ if needed, we can assume that $|\tau(j)|\geq |\Theta (j)|$ for all $j\in R_a$, $y$ remains the fixed point of $\tau$. We define:

\medskip

$\imath )$ 
an alphabet $B$ by  $B =\bigl\{  (j,p); j\in R_a,\;
1\leq p\leq |\Theta (j) | \bigr\} $,
 
$\imath\imath )$
a map $\phi \colon B\to A$  by  $\phi(j,p)= \bigl(\Theta(j)\bigr)_p $ ,

$\imath\imath\imath )$ 
a map $\psi\colon R_a \to B^+ $  by  $\psi(j)=(j,1)(j,2)\ldots (j_,| \Theta (j) | )$.

\medskip        
                                      
The  substitution $\zeta$ on $B$ is defined by: 
$$
\hbox{For $j$ in $R_a$
and $1\leq p\leq | \Theta (j) |$, }\zeta(j,p) =\left\{  \matrix{
\psi\Bigl(\bigl(\tau(j)\bigr)_p\Bigr) & \hbox{ if }&1\leq p< |\Theta (j) |\cr
\psi\Bigl(\bigl(\tau(j)\bigr)_{\textstyle[ |\Theta (j) |, |\tau (j) | ]}\Bigr) & \hbox{ if
}&p= |\Theta (j) |\cr
 }\right.  
$$
The subshift it generates is isomorphic to $(X,T)$, for the details see \cite{DHS} (Section 6).
\hfill $\Box$

\medskip

Hence to know if $\F (X,T)$ or $\F^* (X,T)$ are finite it suffices to compute the substitution $\zeta$ given by the previous proposition and then to use Proposition \ref{subprop}.

\medskip

{\bf An example.}
Let $\sigma$ be the primitive substitution defined on the alphabet $A=\{a,b\}$ by:
$$
\sigma(a)= aba\ ;\ \sigma(b)=baab\ .
$$
Let $(X,T)$ be the non-periodic subshift $\sigma$ generates. We will compute the proper substitution $\zeta$ given by the proof of the previous proposition and then we will deduce that $\F (X,T)$ is infinite and $\F^* (X,T)$ is finite. We will use the notations of the previous proof.

Let $x$ be the one-sided fixed point of $\sigma$ such that $x_{0}=a$:
$$x =  ababaababa \cdots .$$ 

The set of return words to $a$ is $\{ \Theta (1) = ab , \Theta (2) = a \}$. We have:
$$
\sigma ( \Theta (1)) = \sigma(ab)= ababaab= \Theta (1121) \ \ \hbox{\rm and} \ \ \sigma ( \Theta (2) ) = \sigma(a)= aba = \Theta (12)
$$ 
thus 
$$
\tau (1) = 1121 \ \ \hbox{\rm and } \ \ \tau(2)=12.
$$
Let $y\in R_a^{\pNN} = \{ 1,2 \}^{\pNN}$ be the unique one-sided fixed point of $\tau$, it satisfies $\Theta(y)=x$. 

Now we compute the substitution $\zeta : B\rightarrow B^*$. The alphabet $B$ is: 
$$
B=\bigl\{
(1,1),\;(1,2),\;(2,1) \bigr\}
$$ 
The maps $\phi \colon B\to A$ and $\psi\colon \{ 1,2 \}\to B^+$ are given by: $$
\displaylines{
\phi(1,1)=a \ ; \ \phi(1,2)=b \ ; \ \phi(2,1)=a \ ; \cr
\psi(1)=(1,1)(1,2) \ ; \ \psi(2)=(2,1) \ .\cr }
$$ 
We compute now the substitution $\zeta$ on $B$: 
$$
\matrix{ \zeta(1,1)\hfill&=& \psi\bigl(\tau(1)_1\bigr)\hfill &=&
\psi(1)\hfill&=&(1,1)(1,2)\hfill\cr
\zeta(1,2)\hfill&=& \psi\bigl(\tau(1)_{[2,4]}
\bigr)\hfill &=& \psi(121)\hfill&=&(1,1)(1,2)(2,1)(1,1)(1,2)\hfill\cr
\zeta(2,1)\hfill&=&
\psi\bigl(\tau(2)_{[1,2]}\bigr) \hfill&=&
\psi(12)\hfill&=&(1,1)(1,2)(2,1)\hfill\cr
}
$$
Let $M$ be the transposed matrix of $\zeta$, we have 
\[
M        = 
\left[
\begin{array}{ccc}
1       & 1       & 0    \cr 
2       & 2       & 1    \cr
1       & 1       & 1    \cr
\end{array}
\right] , \
e= 
\left(
\begin{array}{ccc}
1   \cr
1   \cr
1   \cr
\end{array}
\right) , \
Me= 
\left(
\begin{array}{cccccc}
2   \cr 
5   \cr
3   \cr
\end{array}
\right) , \    
M^2 e= 
\left(
\begin{array}{cccccc}
7   \cr 
17  \cr
10  \cr
\end{array}
\right) 
\]
We see that these three vectors are linearly independent and that $r=2$. Hence  $Q (X) = X^3 -4X^2 + 2X$. As $\gcd (4,2,0) = 2$ is a prime number, $\F (X,T)$ is infinite and $\F^{*} (X,T)$ is finite.

\bigskip

{\bf Acknowledgments.} The author would like to thank G. Didier and B. Host for helpful discussions.

Fabien Durand

Facult\'e de Math\'ematiques et d'Informatique

Universit\'e de Picardie Jules Verne

33, rue Saint Leu

80039 Amiens Cedex

e-mail : fabien.durand@u-picardie.fr

\end{document}